\newcommand{\M}{{\mathcal M}}    % la varieta' 
\newcommand{\D}{{\mathcal D}}    % la distribuzione
\newcommand{\Hf}{{\vec H}}       % flusso Hamiltoniano
\newcommand{\R}{{I\!\!R}}        % i numeri reali

\documentclass[oneside,draft,11pt]{amsart}
\usepackage{amsmath}        % AmSLaTeX
\usepackage{amsthm}         % aggiunge nuovi ambienti tipo teorema
\usepackage{times}

\title[Normal Extremizers in sub-Riemannian Manifolds]{%
Variational Aspects of the Geodesic Problem in sub-Riemannian Geometry
}
\author{Paolo Piccione \and Daniel V.\ Tausk}
\address{Departamento de Matem\'atica,\hfill\break\indent  
Instituto de Matem\'atica e Estat\'\i stica
\hfill\break\indent Universidade de S\~ao Paulo, Brazil}
\email{piccione@ime.usp.br, tausk@ime.usp.br}
\urladdr{ http://www.ime.usp.br/\~{}piccione,\hfill\break\indent
\phantom{\it URL: }
http://www.ime.usp.br/\~{}tausk}
\thanks{The first author is partially sponsored by CNPq (Brazil) and
the second author is sponsored by FAPESP (S\~ao Paulo, Brazil)}
\date{May 2000}

\theoremstyle{plain}\newtheorem{teo}{Theorem}[section]
\theoremstyle{plain}\newtheorem{lem}[teo]{Lemma}
\theoremstyle{plain}\newtheorem{prop}[teo]{Proposition}
\theoremstyle{plain}\newtheorem{cor}[teo]{Corollary}
\theoremstyle{definition}\newtheorem{defin}[teo]{Definition}
\theoremstyle{remark}\newtheorem{rem}[teo]{Remark}

\begin{document}

\begin{abstract}
We study the local geometry of the space of horizontal
curves with endpoints freely varying in two given submanifolds $\mathcal P$
and $\mathcal Q$ of a manifold $\mathcal M$ endowed with a
distribution $\mathcal D\subset T\M$.
We give a different proof, that holds in a more general context,
of a result by Bismut \cite[Theorem~1.17]{Bis} stating that 
the normal extremizers that are not abnormal are critical points of the
sub-Riemannian  action functional.   We use the Lagrangian
multipliers method in a Hilbert manifold setting,  which leads to a characterization
of the abnormal extremizers (critical points of the endpoint map) as  curves where the
linear constraint fails to be regular. Finally,   we describe a modification of a 
result by Liu and Sussmann \cite{LS} that shows the global distance minimizing
property of sufficiently small portions of normal extremizers between a point and a
submanifold.
\end{abstract}

\maketitle
\begin{section}{Introduction}\label{sec:intro}
%The goal of this paper is to prove some basic results
%concerning the normal extremizers in a semi-Riemannian
%manifold, with the aim of providing the tools to develop
%a variational theory for this kind of geodesics.  
A sub-Riemannian manifold consists of a   smooth $n$-dimensional manifold
$\M$, and a smooth distribution $\D\subset T\M$  on $\M$ of constant 
rank $n-k$, endowed
with  a smoothly varying positive definite metric tensor $g$. 
The length is defined only for {\em horizontal\/}
curves in $\M$, i.e., curves which are everywhere tangent to $\D$.
It was proven  in \cite{LS} that  a horizontal curve which
minimizes length is either a {\em normal extremal\/} or an {\em abnormal
extremal}, where the two possibilities are not mutually exclusive.
This proof is obtained as an application of the {\em Pontryagin maximum
principle\/} of Optimal Control Theory; an alternative proof of
this fact obtained by variational methods is given in this paper
(Corollary~\ref{thm:normalabnormal}).

A normal extremal is defined as a   curve in $\M$ that is a solution
of the sub-Riemannian Hamiltonian $H(p)=\frac12\,g^{-1}(p\vert_\D,
p\vert_\D)$ on $T\M^*$, i.e., a  curve that is
the projection on $\M$ of an integral line of the Hamiltonian
flow $\Hf$. Such curves are automatically horizontal.
An abnormal extremal can be defined as a curve which is the
projection on $\M$ of a non zero {\em characteristic curve\/}
in the annihilator $\D^o\subset T\M^*$; a characteristic curve
is a curve in $\D^o$ which is tangent to the kernel of the restriction
to $\D^o$ of the canonical symplectic form of $T\M^*$.

As in the case of Riemannian geodesics, sufficiently small
segments of a normal extremal is length minimizing (see \cite{LS});
however, ``most'' abnormal extremals do not have any sort
of minimizing property (observe that the definition of abnormal
minimizer does not involve the metric $g$). 

The first example of a length minimizer which is {\em not\/}
a normal extremal was given in \cite{M2}. The goal of this
paper is to discuss the theory of extremals by techniques of
Calculus of Variations and to give the basic instruments to develop
a variational theory (Morse Theory, Ljusternik--Schnirelman
theory) for sub-Riemannian geodesics. The results of this paper
are used in \cite{GGP}, where the authors consider the
problem of existence and multiplicity of geodesics joining a point
and a line 
in a sub-Riemannian manifold $(\M,\D,g)$, with $\mathrm{codim}(\D)=1$.

In \cite[Theorem~1.17]{Bis} it is proven that the normal sub-Riemannian
extremals between two fixed points of a sub-Riemannian manifold
are critical points of the sub-Riemannian action functional.
The proof is presented in the context of the {\em Malliavin calculus},
employed to study some problems connected with the asymptotics of the semi-group
associated with a hypoelliptic diffusion. For this purposes, the author's
proof is restricted to the case that the image of the normal
extremal be contained in an open subset of $\M$ on which the 
distribution $\mathcal D$ is globally generated by $n-k$ smooth vector fields.
In this paper we reprove the result of \cite[Theorem~1.17]{Bis}
under the more general assumptions  that:
\begin{itemize}
\item the vector bundle $\mathcal D$ is not necessarily trivial
around the image of the normal extremizer;
\item the endpoints of the normal extremizers are free to move on two submanifolds
of $\M$.
\end{itemize}
As to the first generalization of the extremizing property of the normal
extremizers, it is interesting to observe that in the proof it is employed
the Lagrangian multipliers technique that uses {\em time-dependent referentials\/}
of $\mathcal D$ defined in a neighborhood of the graph
of any continuous curve in $\M$. The existence of such referentials is obtained 
by techniques of calculus with affine connections, and it is likely that
the method of time-dependent referentials may be applied to other situations
where global geometrical results are to be proven. 
For instance, in \cite{Kishi} the author proves a Morse Index Theorem
for normal extremizers, but in his proof he implicitly
assumes the triviality of the vector bundle $\mathcal D$ in a neighborhood
of the curve. However, the arguments presented could be made more
precise by a systematic use of time-dependent referentials.

Another observation that is worth making about the Lagrangian
multipliers is that, in the functional setup of the method, the constraint
is given by the kernel of a suitable submersion (see formula \eqref{eq:defTheta})
from the set of $H^1$-curves in an open subset of $\M$ taking values
in the Hilbert space of $\R^k$-valued $L^2$-functions. This submersion
is defined using time-dependent referentials of the annihilator $\mathcal D^o$
of $\mathcal D$ in the cotangent bundle $T\M^*$, and the surprising result
is that such map fails to be a submersion precisely at the abnormal
extremizers. We therefore obtain a new variational description of the
abnormal extremizers in a sub-Riemannian manifold.

Finally, it is important to emphasize the role of the endmanifolds
$\mathcal P$ and $\mathcal Q$ in the development of the theory.
An interesting result is that, if either one of the two is everywhere
transversal to $\mathcal D$, then the set of horizontal curves
between $\mathcal P$ and $\mathcal Q$ does not contain singularities
(Proposition~\ref{thm:critpts}); in particular, all the sub-Riemannian
extremizers between $\mathcal P$ and $\mathcal Q$ are normal.
This fact can be used in several circumstances: for instance,
in Corollary~\ref{thm:corballs} we obtain some information
about the geometry of sub-Riemannian balls;  moreover,
it is possible to obtain also some criteria to establish the smoothness
for abnormal extremizers (see Remark~\ref{thm:remsmoothness}).

We outline briefly the contents of each section of this article.

In Section~\ref{sec:diff} we study the local geometry of the
space of horizontal curves joining two fixed points $q_0$ and $q_1$
of $\M$ by two different techniques. On one hand, this space can
be described as the set of curves $\gamma$ joining $q_0$ and $q_1$
satisfying $\theta_i(\dot\gamma)=0$, where $\theta_1,\ldots,\theta_k$
is a local time-dependent referential for the annihilator
$\D^o$ of $\D$. On the other hand, the same space can be obtained
as the inverse image of $q_1$  by the {\em endpoint mapping\/}
restricted to the set of horizontal curves emanating from $q_0$.
We show that these two constraints have the same regular points;
such curves are called {\em regular\/} and a suitable neighborhood
of them in the space of horizontal curves joining $q_0$ and $q_1$
has the structure of an infinite dimensional Hilbert manifold.

In Section~\ref{sec:normalgeo} we define the normal extremals,
also called {\em normal geodesics}, in a sub-Riemannian manifold,
using the Hamiltonian setup.

In Section~\ref{sec:abnormal} we study the image of the differential
of the endpoint mapping; to this aim we introduce an atlas on the
space of horizontal curves starting at $q_0$.

Finally, in Section~\ref{sec:critical} we prove that a regular
curve is a critical point of the sub-Riemannian action functional
if and only if it is a normal geodesic. We also study the case
of curves with endpoints varying in two submanifolds of $\M$.
If we consider the space of horizontal curves joining the
submanifolds $\mathcal P$ and $\mathcal Q$, then, provided that
either $\mathcal P$ or $\mathcal Q$ is transversal to $\D$, this
set is always a Hilbert manifold. Moreover, the critical points
of the sub-Riemannian action functional in this space are those
normal geodesics between $\mathcal P$ and $\mathcal Q$ whose
Hamiltonian lift annihilates the tangent spaces of $\mathcal P$
and $\mathcal Q$ at its endpoints.

To conclude the paper, we present two short appendices. 
In Appendix~\ref{sec:affine} we prove that every
horizontal curve can be obtained as the reparameterization of
an affinely parameterized horizontal curve. In Appendix~\ref{sec:localminimality}
we adapt a proof of local optimality of
normal geodesics due to Liu and Sussmann \cite[Appendix~C]{LS}
to prove that sufficiently small portions of normal geodesics 
are length minimizers between an initial submanifold and a point.

\end{section}

\begin{section}{The differentiable structure of the space of\\ horizontal
curves}\label{sec:diff}
We give a couple of preliminary results needed to the study
of the geometry of the set of horizontal paths in a sub-Riemannian 
manifold.
The main reference for the geometry of infinite dimensional
manifolds is \cite{L}; for the basics of Riemannian geometry
we refer to \cite{dC}.

Recall that a smooth map $f:M\mapsto N$
between Hilbert manifolds is a {\em submersion\/} at $x\in M$ if
the differential  $\mathrm df(x):T_xM\mapsto T_{f(x)}N$ is surjective;
$f$ is a submersion if it is a submersion at every $x\in M$.

\begin{lem}\label{thm:subfg}
Let $M$, $M_1$ and $M_2$ be Hilbert manifolds and let $f:M\mapsto M_1$,
$g:M\mapsto M_2$ be submersions. Let $p_1\in M_1$, $p_2\in M_2$ and choose
$x\in f^{-1}(p_1)\cap g^{-1}(p_2)$. Then, $f\vert_{g^{-1}(p_2)}$
is a submersion at $x$ if and only if $g\vert_{f^{-1}(p_1)}$ is a submersion
at $x$.
\end{lem}
\begin{proof}
We need to show that $\mathrm df(x)\vert_{\mathrm{Ker}(\mathrm dg(x))}$
is surjective onto $T_{f(x)}M_1 $ if and only if 
$\mathrm dg(x)\vert_{\mathrm{Ker}(\mathrm df(x))}$
is surjective onto $T_{g(x)}M_2 $. This follows from a general fact:
if $T:V\mapsto V_1$ and $S:V\mapsto V_2$ are surjective
linear maps between
vector spaces, then $T\vert_{\mathrm{Ker}(S)}$ is surjective if and only if
$\mathrm{Ker}(T)+\mathrm{Ker}(S)=V$. Clearly, this relation is symmetric
in $S$ and $T$, and we obtain the thesis.
\end{proof}

We give one more introductory result concerning the existence
of time-depen\-dent local referentials for vector bundles defined in 
a neighborhood of a given curve. We need the following definition:
\begin{defin}\label{thm:defnormalradius}
Let $(\M,\overline g)$ be a Riemannian manifold and $x\in \M$. A positive
number $r\in\R^+$ is said to be a {\em normal radius\/}
for $x$ if $\exp_x:{\mathrm B}_r(0)\mapsto B_r(x)$ is a diffeomorphism,
where $\exp$ is the exponential map of $(\M,\overline g)$, ${\mathrm B}_r(0)$
is the open ball of radius $r$ around $0\in T_x\M$ and
$B_r(x)$ is the open ball of radius $r$ around $x\in\M$.
We say that $r$ is {\em totally normal\/} for $x$ if $r$ is a normal
radius for all $y\in B_r(x)$.
\end{defin}
By a simple argument in Riemannian geometry, it is easy to see
that if $K\subset\M$ is a compact subset, then there exists $r>0$
which is totally normal for all $x\in K$.
\smallskip

Given an  vector bundle $\pi:\xi\mapsto\M$ of rank $k$
over a manifold $\M$, a {\em time-dependent local referential\/} of
$\xi$ is a family of smooth maps $X_i:U\mapsto\xi$, $i=1,\ldots,k$,
defined on an open subset $U\subseteq\R\times\M$ such that 
$\{X_i(t,x)\}_{i=1}^k$ is a basis of the fiber $\xi_x$ for all
$(t,x)\in U$.
\begin{lem}\label{thm:basis}
Let $\M$ be a finite dimensional manifold, let $\pi:\xi\mapsto\M$ be
a vector bundle over $\M$ and let $\gamma:[a,b]\mapsto\M$ be a continuous
curve. Then, there exists an open subset $U\subseteq  \R\times\M$ containing
the graph of $\gamma$ and a smooth time-dependent local referential
of $\xi$ defined in $U$.
\end{lem}
\begin{proof}
We first consider the case that $\gamma$ is a smooth curve.
Let us choose an arbitrary connection in $\xi$, an arbitrary Riemannian
metric $\overline g$ on $\M$ and   a smooth
extension $\gamma:[a-\varepsilon,b+\varepsilon]\mapsto\M$ of
$\gamma$, with $\varepsilon>0$. 
Since the image of $\gamma$ is compact in $\M$, there
exists $r>0$ which is a normal radius for all $\gamma(t)$, $t\in[a-\varepsilon,
b+\varepsilon]$. We define $U$ to be the open set:
\[U=\Big\{(t,x)\in\R\times\M:t\in\,]a-\varepsilon,b+\varepsilon[,\ 
x\in B_r(\gamma(t))\Big\}.\]
Let now $\overline X_1,\ldots,\overline X_k$ be a referential
of $\xi$ along $\gamma$; for instance, this referential can
be chosen by parallel transport along $\gamma$ relative to the
connection on $\xi$. Finally, we obtain a time-dependent
local referential for $\xi$ in $U$ by setting, for
$(t,x)\in U$ and for $i=1,\ldots,k$, $X_i(t,x)$ equal to
the parallel transport (relative to the connection of $\xi$)
of $\overline X_i(t)$ along the radial geodesic joining $\gamma(t)$ and
$x$.

The general case of a continuous curve is easily obtained by
a density argument. For, let $\gamma:[a,b]\mapsto\M$
be continuous and let $r>0$ be a totally normal radius for
$\gamma(t)$, for all $t\in[a,b]$. Let $\gamma_1:[a,b]\mapsto\M$
be any smooth curve such that $\mathrm{dist}(\gamma(t),\gamma_1(t))<r$
for all $t$, where $\mathrm{dist}$ is the distance induced
by the Riemannian metric $\overline g$ on $\M$. Then, if we repeat the
above proof for the curve $\gamma_1$, the open set $U$ thus obtained
will contain the graph of $\gamma$, and we are done.
\end{proof}
Let us now consider a sub-Riemannian manifold, that is a triple $(\M,\D,g)$
where $\M$ is a smooth $n$-dimensional manifold, $\D$ is a smooth distribution
in $\M$ of codimension $k$ and $g$ is smoothly varying positive inner
product on $\D$.

A curve $\gamma:[a,b]\mapsto\M$ is said
to be {\em $\D$-horizontal}, or simply {\em horizontal}, if
it is absolutely continuous and if
$\dot\gamma(t)\in\D$ for almost all $t\in[a,b]$. As we did in the proof of
Lemma~\ref{thm:basis}, we will use sometimes auxiliary structures
on $\M$, which are chosen (in a non canonical way) once for
all. We therefore assume that $\overline g$ is a given Riemannian
metric tensor on $\M$ such that $\overline g\vert_\D=g$, that
$\D_1$ is a $k$-dimensional distribution in $\M$ which is
complementary to $\D$ (for instance, $\D_1$ is the $\overline g$-orthogonal
distribution to $\D$), and we also assume that $\nabla$ is a linear
connection in $T\M$ which is {\em adapted\/} to the decomposition
$\D\oplus\D_1$, i.e., the covariant derivative of  vector fields
in $\D$ (resp., in $\D_1$) belongs to $\D$ (resp., to $\D_1$). 
For the construction of these objects, one can consider an arbitrary
Riemannian metric $\widetilde g$ on $\M$. Then, one defines $\D_1$ 
as the $\widetilde g$-orthogonal complement of $\D$ and $\overline
g\vert_{\D_1}=
\widetilde g\vert_{\D_1}$; for the connection $\nabla$, it suffices to choose 
any pair of connections $\nabla_0$ and $\nabla_1$ respectively
on the vector bundles $\D$ and $\D_1$ and then one sets
$\nabla=\nabla_0\oplus\nabla_1$. Observe that the connection $\nabla$
constructed in this way is {\em not\/} torsion free; we denote by $\tau$
the torsion of $\nabla$:
\[\tau(X,Y)=\nabla_XY-\nabla_YX-[X,Y].\]

Using Lemma~\ref{thm:basis}, we describe $\D$ locally as the kernel
of a time-dependent $\R^k$-valued $1$-form:
\begin{prop}\label{thm:kform}
Let $\gamma:[a,b]\mapsto\M$ be a continuous curve. Then, there
exists an open subset $U\subseteq\R\times\M$ containing the graph
of $\gamma$ and a smooth time-dependent $\R^k$-valued
$1$-form $\theta$ defined in $U$,
with
$\theta_{(t,x)}:T_x\M\mapsto\R^k$  a surjective linear map and
$\D_x=\mathrm{Ker}(\theta_{(t,x)})$ for all $(t,x)\in U$.
\end{prop}
\begin{proof}
Let $\xi$ be the subbundle of the cotangent bundle $T\M^*$
given by the {\em annihilator\/} $\D^o$ of $\D$. Apply Lemma~\ref{thm:basis}
to $\xi$ and set $\theta=(\theta_1,\dots, \theta_k)$, where 
$\{\theta_i\}_{i=1}^k$ is a time-dependent 
local referential of $\xi$ defined in an open 
neighborhood of the graph of $\gamma$.
\end{proof}
Observe that, since $\D_1$ is complementary to $\D$, for
all $(t,x)\in U$ the map \[\theta_{(t,x)}:\D_1\longmapsto\R^k\]
is an isomorphism.

Let us now consider the following spaces of curves in $\M$.

We denote by $L^2([a,b],\R^m)$ the Hilbert space of Lebesgue
square integrable $\R^m$-valued maps on $[a,b]$ and by $H^1([a,b],\R^m)$
the Sobolev space of all absolutely continuous maps $x:[a,b]\mapsto\R^m$
with derivative in $L^2([a,b],\R^m)$. Finally, we denote
by $H^1([a,b],\M)$ the set of curves $x:[a,b]\mapsto\M$ such that 
for any local chart $(U,\phi)$ on $\M$, with $\phi:U\subset\M\mapsto\R^n$,
and for any closed interval $I\subset x^{-1}(U)$, the map
$\phi\circ( x\vert_I)$ is in $H^1(I,\R^m)$. It is well known
that $H^1([a,b],\M)$ is an infinite dimensional smooth manifold modeled
on the Hilbert space $H^1([a,b],\R^n)$ (see for instance \cite{PT1}
for a recent reference on these issues). 

For all pairs of points $q_0,q_1\in\M$, we define the following sets
of curves in $\M$:
\begin{equation}\label{eq:defspaces}
\begin{split}
&
H^1_{q_0}([a,b],\M)=\Big\{x\in H^1([a,b],\M):x(a)=q_0\Big\};\\
& H^1_{q_0,q_1}([a,b],\M)=\Big\{x\in H^1([a,b],\M):x(a)=q_0,\ x(b)=q_1\Big\};\\
& H^1([a,b],\D,\M)=\Big\{x\in H^1([a,b],\M):\dot x(t)\in\D\ \text{a.e.\
on $[a,b]$}\Big\};
\\
& H^1_{q_0}([a,b],\D,\M)=H^1([a,b],\D,\M)\cap H^1_{q_0}([a,b],\M);\\
&H^1_{q_0,q_1}([a,b],\D,\M)=H^1([a,b],\D,\M)\cap H^1_{q_0,q_1}([a,b],\M).
\end{split}
\end{equation}
We prove that  the sets $H^1_{q_0}([a,b],\M)$,
$H^1_{q_0,q_1}([a,b];\M)$, $H^1([a,b],\D,\M)$ and
$H^1_{q_0}([a,b],\D,\M)$,
are smooth submanifolds of $H^1([a,b],\M)$
for all $q_0,q_1\in\M$.
However, in general, the space $H^1_{q_0,q_1}([a,b],\D,\M)$, 
consisting of horizontal curves
joining the two fixed points $q_0$ and $q_1$, is {\em not\/}
a submanifold of $H^1_{q_0,q_1}([a,b];\M)$, and this fact is precisely
the origin of difficulties when one tries to develop a variational
theory for sub-Rieman\-nian geodesics. 

\noindent\ \ 
In  order to see that $H^1_{q_0}([a,b],\M)$ and $H^1_{q_0,q_1}([a,b],\M)$ are
submanifolds of $H^1([a,b],\M)$, simply observe that the map \[\mathcal
E_{a,b}:\gamma\mapsto(\gamma(a),
\gamma(b))\] is a submersion of $H^1([a,b],\M)$ into
$\M\times\M$. 

Then,
$H^1_{q_0}([a,b],\M)=\mathcal E^{-1}_{a,b}(\{q_0\}\times\M)$ and
$H^1_{q_0,q_1}([a,b],\M)=\mathcal E^{-1}_{a,b}(q_0,q_1)$ are smooth submanifolds
of $H^1([a,b],\M)$.
\smallskip

As to the regularity of $H^1_{q_0}([a,b],\D,\M)$, we will now show that
this set can be covered by a family of open subset $\{\mathcal U_\alpha\}$
of $H^1_{q_0}([a,b],\M)$ such that each intersection
$H^1_{q_0}([a,b],\D,\M)\cap \mathcal U_\alpha$
is the inverse image of a submersion of $\mathcal U_\alpha$ in 
the Hilbert space $L^2([a,b],\R^k)$. The regularity of
$H^1([a,b],\D,\M)$ will follow by a similar argument.

To this aim, 
let $\gamma_0$ be a fixed curve in $H^1_{q_0}([a,b],\M)$ and let
$U_{\gamma_0}\subset\R\times\M$ be an open set containing the graph of
$\gamma_0$ and that is the domain of the map $\theta$
of Proposition~\ref{thm:kform}. Denote by $H^1_{q_0}([a,b],\M,U_{\gamma_0})$
the open subset of $H^1_{q_0}([a,b],\M)$ consisting of those curves
whose graphs is contained in $U_{\gamma_0}$:
\begin{equation}\label{eq:defUgamma0}
H^1_{q_0}([a,b],\M,U_{\gamma_0})=\Big\{\gamma\in
H^1_{q_0}([a,b],\M):(t,\gamma(t))\in U_{\gamma_0},\ \text{for all}\
t\in[a,b]\Big\}.
\end{equation}
Let $\Theta:H^1_{q_0}([a,b],\M,U_{\gamma_0})\mapsto L^2([a,b],\R^k)$ be the
smooth map defined by:
\begin{equation}\label{eq:defTheta}
\Theta(\gamma)(t)=\theta_{(t,\gamma(t))}(\dot\gamma(t)).
\end{equation}
Clearly, $H^1_{q_0}([a,b],\M,U_{\gamma_0})\cap
H^1_{q_0}([a,b],\D,\M)=\Theta^{-1}(0)$.
\begin{prop}\label{thm:Thetasub}
$\Theta$ is a submersion.
\end{prop}
\begin{proof}
Clearly $\Theta$ is smooth because $\theta$ is smooth.
To compute the differential of $\Theta$ we use the connection $\nabla$
adapted to the decomposition $T\M=\D\oplus\D_1$ introduced above.
Let $\gamma\in H^1_{q_0}([a,b],\M,U_{\gamma_0})$ be fixed and let $V\in T_\gamma
H^1_{q_0}([a,b],\M)$, i.e., $V$ is a vector field of class $H^1$ along
$\gamma$ with $V(a)=0$. We write $V=V_\D+V_{\D_1}$ with
$V_\D(t)\in\D$ and $V_{\D_1}(t)\in\D_1$ for all $t$; using the properties of
$\nabla$ we compute easily:
\begin{equation}\label{eq:dTheta}
\mathrm d\Theta(\gamma)[V](t)=\!\left[\nabla_V\theta\right]_{(t,\gamma(t))}
(\dot\gamma(t))+\theta_{(t,\gamma(t))}(\nabla_{\dot\gamma(t)}V)
+\theta_{(t,\gamma(t))}(\tau(V(t),\dot\gamma(t))),
\end{equation}
where $\nabla_V\theta$ is the covariant derivative of $\theta_{(t,\cdot)}$.

Let now $f\in L^2([a,b],\R^k)$ be fixed; for the surjectivity of
$\mathrm d\Theta(\gamma)$ we want to solve the
equation in $V$: $\mathrm d\Theta(\gamma)[V]=f$. To this aim,
we choose $V_{\D_0}=0$, and we get:
\begin{equation}\label{eq:eqinV1}
\theta_{(t,\gamma(t))}(\nabla_{\dot\gamma(t)}V_{\D_1})+
\left[\nabla_{V_{\D_1}}\theta\right]_{(t,\gamma(t))}(\dot\gamma(t))+
\theta_{(t,\gamma(t))}(\tau(V_{\D_1}(t),\dot\gamma(t)))=f.
\end{equation}
Since $\theta_{(t,\gamma(t))}:(\D_1)_{\gamma(t)}\mapsto\R^k$ is
an isomorphism, \eqref{eq:eqinV1} is equivalent to
a first order linear differential equation in $V_{\D_1}$, that
admits a unique solution satisfying $V_{\D_1}(a)=0$.
Observe that since $\gamma\in H^1([a,b],\M)$, by \eqref{eq:eqinV1}
we get that $V$ is also of class $H^1$, and we are done.
\end{proof}
\begin{cor}\label{thm:H1pDsubmanifold}
$H^1([a,b],\D,\M)$ and
$H^1_{q_0}([a,b],\D,\M)$ are smooth submanifolds of $H^1([a,b],\M)$.\qed
\end{cor}
We now consider the {\em endpoint mapping\/} $\mathrm{end}:H^1_{q_0}([a,b],\M)
\mapsto\M$ given by:
\[\mathrm{end}(\gamma)=\gamma(b).\]
It is easy to see that $\mathrm{end}$ is a submersion, hence we
have the following:
\begin{cor}\label{thm:Thetasubendsub}
Let $\gamma_0\in H^1_{q_0}([a,b],\M)$ be fixed and let
$H^1_{q_0}([a,b],\M,U_{\gamma_0})$,
$\Theta$ be defined as in \eqref{eq:defUgamma0} and \eqref{eq:defTheta}.

Then, for all $\gamma\in\Theta^{-1}(0)\cap\mathrm{end}^{-1}(q_1)=
H^1_{q_0,q_1}([a,b],\D,\M)$, the restriction
$\Theta\vert_{H^1_{q_0}([a,b],\M,U_{\gamma_0})\cap H^1_{q_0,q_1}([a,b],\M) }$
is a submersion if and only if the restriction
$\mathrm{end}\vert_{H^1_{q_0}([a,b],\D,\M)}$ is a submersion.
\end{cor}
\begin{proof}
It follows immediately from Lemma~\ref{thm:subfg} 
and Proposition~\ref{thm:Thetasub}.
\end{proof}
\begin{defin}\label{thm:defregular}
A curve $\gamma\in H^1_{q_0,q_1}([a,b],\D,\M)$ is said to be
{\em regular\/} if the restriction $\mathrm{end}\vert_{H^1_{q_0}([a,b],\D,\M)}$ 
is a submersion at $\gamma$. If $\gamma$ is not regular, then it  is called
an {\em abnormal extremal}.
\end{defin}
Observe that the notion of abnormal extremality is not related
to any sort of extremality with respect to the length or the
action functional, but rather to lack of regularity in
the geometry of the space of horizontal paths.
The smoothness of length minimizing abnormal extremals is an open question.

\end{section}

\begin{section}{Normal Geodesics}\label{sec:normalgeo}

In order to define the normal geodesics in a sub-Riemannian
manifold we introduce a Hamiltonian setup in $T\M^*$ as follows. 

Let us consider the cotangent bundle $T\M^*$ endowed with its canonical
symplectic form $\omega$. Recall that $\omega$ is defined by $\omega=
-\mathrm d\vartheta$, $\vartheta$ being the canonical $1$-form on $T\M^*$
given by $\vartheta_p(\rho)=p(\mathrm d\pi_p(\rho))$, where
$\pi:T\M^*\mapsto\M$ is the projection, $p\in T\M^*$ and $\rho\in T_pT\M^*$.
Let $H:T\M^*\mapsto\R$ be a smooth function; we call such a function
a {\em Hamiltonian\/} in $(T\M^*,\omega)$. The {\em Hamiltonian
vector field\/} of $H$ is the smooth vector field on $T\M^*$ denoted
by $\Hf$ and defined by the relation $\mathrm dH(p)=\omega(\Hf(p),\cdot)$;
the integral curves of $\Hf$ are called the {\em solutions\/} of the 
Hamiltonian $H$. With a slight abuse of terminology, we will say that a smooth
curve $\gamma:[a,b]\mapsto\M$ is a solution of the Hamiltonian $H$
if it admits a {\em lift\/} $\Gamma:[a,b]\mapsto T\M^*$ that is
a solution of $H$.

More in general, one can consider {\em time-dependent\/} Hamiltonian
functions on $T\M^*$, which are smooth maps defined
on an open subset $U$ of $\R\times T\M^*$. In this case,
the Hamiltonian flow $\Hf$ is a time-dependent vector field in $T\M^*$,
and its integral curves in $T\M^*$ are again called the solutions of
the Hamiltonian $H$.

A {\em symplectic chart\/} in $T\M^*$ is a local chart  taking
values in $\R^n\oplus{\R^n}^*$ whose differential at each point is
a symplectomorphism from the tangent space $T_p(T\M^*)$ to
$\R^n\oplus{\R^n}^*$ endowed with the canonical symplectic structure.
 Given a chart $q=(q_1,\ldots, q_n)$ in $\M$, 
we get a symplectic chart $(q,p)$ on
$T\M^*$ where $p=(p_1,\ldots,p_n)$ is defined by
$p_i(\alpha)=\alpha\left(\frac\partial{\partial q_i}\right)$. 
We denote by
$\left\{\frac\partial{\partial q_i},\frac\partial{\partial p_j}
\right\}$, $i,j=1,\ldots,n$, the corresponding local referential
for $T(T\M^*)$, and by $\{\mathrm dq_i,\mathrm dp_j\}$ the local referential
of $T(T\M^*)^*$. We have:
\[\omega=\sum_{i=1}^n\mathrm dq_i\wedge\mathrm dp_i,\quad
\Hf=\sum_{i=1}^n\left(\frac{\partial H}{\partial p_i}\,
\frac\partial{\partial q_i}-
\frac{\partial H}{\partial q_i}\,\frac\partial{\partial p_i}\right).\]
In the symplectic chart $(q,p)$, a solution $\Gamma(t)=(q(t),p(t))$
of the Hamiltonian $H$ is the solution of the Hamilton  equations:
\begin{equation}\label{eq:HJeqns}
\left\{\begin{array}{l}\displaystyle \frac{\mathrm dq}{\mathrm dt}=\phantom{-}\frac{
\partial H}{\partial p},\\  \\ \displaystyle
\frac{\mathrm dp}{\mathrm dt}=-\frac{\partial H}{\partial q}.\end{array}\right.
\end{equation}

\begin{defin}\label{thm:defnormgeo}
A {\em normal geodesic\/} in the sub-Riemannian manifold
 $(\M,\D,g)$ is a curve $\gamma:[a,b]\mapsto\M$ that admits a lift
$\Gamma:[a,b]\mapsto T\M^*$ which is a solution of the sub-Riemannian
Hamiltonian $H:T\M^*\mapsto\R$ given by:
\begin{equation}\label{eq:defHsR}
H(p)=\frac12\,g^{-1}(p\vert_\D,p\vert_\D),
\end{equation}
where $g^{-1}$ is the induced inner product in $\D^*$. In this case,
we say that $\Gamma$ is a {\em Hamiltonian lift\/} of $\gamma$.
\end{defin}
The Hamilton  equations for
 the sub-Riemannian Hamiltonian
\eqref{eq:defHsR} will be computed explicitly in Section~\ref{sec:critical}
(formula~\eqref{eq:newHJH}). It will be seen that the first of the two equations
means that the solutions in
$\M$ are horizontal curves and that $\Gamma\vert_{\D}=g(\dot\gamma,\cdot)$
(see remark~\ref{thm:remlift}).

We remark that a normal geodesic need not be regular in the sense
of Definition~\ref{thm:defregular}, 
hence there are geodesics that are at the same
time normal and abnormal.
Observe also that, in general, a normal geodesic $\gamma$ may 
admit more than one Hamiltonian lift $\Gamma$. 
This phenomenon occurs precisely
when $\gamma$ is at the same time a normal geodesic and
an abnormal extremizer.
\end{section}

\begin{section}{Abnormal extremals and the endpoint mapping}\label{sec:abnormal}
In this section we give necessary and sufficient conditions 
for a curve to be an abnormal extremal in terms of the symplectic structure of
the cotangent bundle $T\M^*$. We describe a coordinate system in the Hilbert
manifold $H^1_{q_0}([a,b],\M)$ which is compatible with the submanifold
$H^1_{q_0}([a,b],\D,\M)$. This will provide an explicit description of the
tangent space $T_\gamma H^1_{q_0}([a,b],\D,\M)$ which will allow us to compute
the image of the differential of the restriction of the endpoint mapping to
$H^1_{q_0}([a,b],\D,\M)$.

Let $\M$ be a manifold endowed with a distribution $\D$, with
$\mathrm{dim}(\M)=n$ and $\mathrm{codim}(\D)=k$. The
sub-Riemannian metric will be irrelevant in the theory of this section. Let
$U\subset\R\times\M$ be an open set and let $X_1,\ldots,X_n$ be a
time-dependent referential of $T\M$ defined in $U$. We say that such
referential is  {\em adapted\/} to the distribution $\D$ if
$X_1,\ldots,X_{n-k}$ form a referential for $\D$.

It follows easily from Lemma~\ref{thm:basis} that, given a continuous 
curve $\gamma:[a,b]\mapsto\M$, there exists an open set $U\subset\R\times\M$
containing the graph of $\gamma$ and a referential of $T\M$ defined in $U$
which is adapted to $\D$. Namely, one chooses a vector subbundle $\D_1\subset
T\M$ such that $T\M=\D\oplus\D_1$ and then apply Lemma~\ref{thm:basis} to both
$\D$ and $\D_1$.

Given a time-dependent referential of $T\M$ defined in an open set 
$U\subset\R\times\M$, we are going to associate to it a map
\[\mathcal B:H^1([a,b],\M,U)\longmapsto L^2([a,b],\R^n),\]
where $H^1([a,b],\M,U)$ denotes the open set in $H^1([a,b],\M)$ consisting 
of curves whose graph is contained in $U$. We define
$\mathcal B$ by:
\begin{equation}\label{eq:defB}
\mathcal B(\gamma)=h,
\end{equation}
where $h=(h_1,\ldots,h_n)$ is given by \begin{equation}\label{eq:defh}
\dot\gamma(t)=\sum_{i=1}^nh_i(t)X_i(t,\gamma(t)),
\end{equation}
for almost all $t\in[a,b]$. The map $\mathcal B$ is smooth. It's differential is computed in the following:

\begin{lem}\label{thm:diffB}
Let $\gamma\in H^1([a,b],\M,U)$ and $v$ be an $H^1$ vector field along $\gamma$. Set $h=\mathcal B(\gamma)$, $z=\mathrm d\mathcal B_\gamma(v)$. We define a time-dependent vector field in $U$ by
\begin{equation}\label{eq:defX}
X(t,x)=\sum_{i=1}^nh_i(t)X_i(t,x),\quad (t,x)\in U
\end{equation}
and a vector field $w$ along $\gamma$ by \begin{equation}\label{eq:defw}
w(t)=\sum_{i=1}^nz_i(t)X_i(t,\gamma(t)).
\end{equation}
Given a chart $(q_1,\ldots,q_n)$ defined in an open set $V\subset\M$, denote by $\tilde v(t)$, $\tilde X(t,q)$ and $\tilde w(t)$ the representation in coordinates of $v$, $X$ and $w$ respectively. Then, the following relation holds:
\begin{equation}\label{eq:difeqvz}
\frac{\mathrm d}{\mathrm dt}\tilde v(t)=\frac{\partial\tilde X}{\partial q}(t,\gamma(t))\tilde v(t)+\tilde w(t),
\end{equation}
for all $t\in[a,b]$ such that $\gamma(t)\in V$.
\end{lem}
\begin{proof}
Simply consider a variation of $\gamma$ with variational vector field $v$ and differentiate relation \eqref{eq:defh} with respect to the variation parameter, using the local chart.
\end{proof}
\begin{cor}
The restriction of the map $\mathcal B$ to the set
\[H^1_{q_0}([a,b],\M,U)=H^1_{q_0}([a,b],\M)\cap H^1([a,b],\M,U)\]
is a local chart, taking values in an open subset of $L^2([a,b],\R^n)$.
\end{cor}
\begin{proof}
For $\gamma\in H^1_{q_0}([a,b],\M)$ the tangent space $T_\gamma H^1_{q_0}([a,b],\M)$ consists of those $H^1$ vector fields $v$ along $\gamma$ such that $v(a)=0$. For a fixed $z\in L^2([a,b],\R^n)$, formula \eqref{eq:difeqvz} is a first order linear differential equation for $\tilde v$; Lemma~\ref{thm:diffB} and standard results of existence and uniqueness of solutions of linear differential equations imply that the differential of $\mathcal B$ at any $\gamma\in H^1_{q_0}([a,b],\M,U)$ maps 
the tangent space $T_\gamma H^1_{q_0}([a,b],\M)$ isomorphically onto
$L^2([a,b],\R^n)$. It follows from the inverse function theorem that $\mathcal
B$ is a local diffeomorphism in $H^1_{q_0}([a,b],\M,U)$. Finally, by standard
results on uniqueness of solutions of differential equations, we see that the
restriction of $\mathcal B$ to $H^1_{q_0}([a,b],\M,U)$ is injective.
\end{proof}

If the referential $X_1,\ldots,X_n$ defining $\mathcal B$ is adapted to $\D$,
then  a curve $\gamma$ in $H^1_{q_0}([a,b],\M,U)$ is horizontal if and only if
$\mathcal B(\gamma)=h$ satisfies $h_{n-k+1}=\ldots=h_n=0$. This means that
$\mathcal B$ is a {\em submanifold chart\/} for $H^1_{q_0}([a,b],\D,\M)$. This
observation will provide a good description of the tangent space $T_\gamma
H^1_{q_0}([a,b],\D,\M)$.

Let $\gamma\in H^1_{q_0}([a,b],\M,U)$ and set $h=\mathcal B(\gamma)$. Define a time-dependent vector field $X$ in $U$ as in \eqref{eq:defX}. By Lemma~\ref{thm:diffB}, the kernel $\mathrm{Ker}\,\mathrm d\mathcal B_\gamma$ is the vector subspace of $T_\gamma H^1([a,b],\M)$ consisting of those $v$ whose representation in coordinates $\tilde v$ satisfy the homogeneous part of the linear differential equation \eqref{eq:difeqvz}, namely:
\begin{equation}\label{eq:fundamentaldif}
\frac{\mathrm d}{\mathrm dt}\tilde v(t)=\frac{\partial\tilde X}{\partial q}(t,\gamma(t))\tilde v(t).
\end{equation}
By the uniqueness of the solution of a Cauchy problem,
it follows that, for all
$t\in[a,b]$, the evaluation map
\[\mathrm{Ker}\,\mathrm d\mathcal B_\gamma\ni v\mapsto v(t)\in T_{\gamma(t)}\M\]
is an isomorphism. Therefore, for every $t\in[a,b]$ we can define a linear isomorphism $\Phi_t:T_{\gamma(a)}\M\mapsto T_{\gamma(t)}\M$ by:
\begin{equation}\label{eq:defPhi}
\Phi_t(v(a))=v(t),\quad v\in\mathrm{Ker}\,\mathrm d\mathcal B_\gamma.
\end{equation}
Using the maps $\Phi_t$ we can give a coordinate free description of the differential of $\mathcal B$, based on the ``method of variation of constants'' for solving non homogeneous linear differential equations.

\begin{lem}\label{thm:diffBcoordfree}
Let $\gamma\in H^1_{q_0}([a,b],\M,U)$ and $v\in T_\gamma H^1_{q_0}([a,b],\M)$. Set $h=\mathcal B(\gamma)$ and $z=\mathrm d\mathcal B_\gamma(v)$. Define the objects $X$, $w$ and $\Phi_t$ as in \eqref{eq:defX}, \eqref{eq:defw} and \eqref{eq:defPhi} respectively. Then, the following equality holds:
\begin{equation}\label{eq:diffBcoordfree}
v(t)=\Phi_t\int_a^t\Phi_s^{-1}w(s)\mathrm ds.
\end{equation}
\end{lem}
\begin{proof}
The right side of \eqref{eq:diffBcoordfree} vanishes at $t=a$, therefore, to conclude the proof, one only has to show that its representation in local coordinates satisfies the differential equation \eqref{eq:difeqvz}. This follows by direct computation, observing that the representation in local coordinates of the maps $\Phi_t$ is a solution of the homogeneous linear differential equation \eqref{eq:fundamentaldif}.
\end{proof}
\begin{cor}\label{thm:tangspace}
Suppose that the referential $X_1,\ldots,X_n$
defining $\mathcal B$ is adapted to
$\D$.  Let $\gamma$  be an horizontal curve in $ H^1_{q_0}([a,b],\M,U)$. Then,
the tangent space $T_\gamma H^1_{q_0}([a,b],\D,\M)$ consists of all vector
fields $v$ of the form
\eqref{eq:diffBcoordfree}, where $w$ runs over all $L^2$ horizontal vector
fields along $\gamma$.
\end{cor}
\begin{proof}
Follows directly from Lemma~\ref{thm:diffBcoordfree}, observing that $\mathcal B$ is a submanifold chart for $H^1_{q_0}([a,b],\D,\M)$, as remarked earlier.
\end{proof}

We now relate the differential of the endpoint map with the symplectic structure of $T\M^*$. We denote by $\D^o\subset T\M^*$ the annihilator of $\D$. The restriction $\omega\vert_{\D^o}$ of the canonical symplectic form of $T\M^*$ to $\D^o$ is in general no longer nondegenerate and its kernel $\mathrm{Ker}(\omega\vert_{\D^o})(p)$ at a point $p\in\D^o$ may be non zero. We say that an absolutely continuous curve $\eta:[a,b]\mapsto\D^o$ is a {\em characteristic curve\/} for $\D$ if
\[\dot\eta(t)\in\mathrm{Ker}(\omega\vert_{\D^o})(\eta(t)),\]
for almost all $t\in[a,b]$.

We take a closer look at the kernel of $\omega\vert_{\D^o}$. Let $Y$ be a horizontal vector field in an open subset of $\M$. We associate to it a
Hamiltonian function $H_Y$ defined by
\[H_Y(p)=p(Y(x)),\]
where $x=\pi(p)$. We can now compute the $\omega$-orthogonal complement of $T_p\D^o$ in $T_pT\M^*$. Recall that $\Hf_Y$ denotes the corresponding
Hamiltonian vector field in $T\M^*$.

\begin{lem}\label{thm:orthcompl}
Let $p\in T\M^*$ and set $x=\pi(p)$. The $\omega$-orthogonal complement of $T_p\D^o$ in $T_pT\M^*$ is mapped isomorphically by $\mathrm d\pi_p$ onto $\D_x$. Moreover, if $Y$ is a horizontal vector field defined in an open neighborhood of $x$ in
$\M$, then $\Hf_Y(p)$ is the only vector in the $\omega$-orthogonal complement
of $T_p\D^o$ which is mapped by $\mathrm d\pi_p$ into $Y(x)$.
\end{lem}
\begin{proof}
The function $H_Y$ vanishes on $\D^o$ and therefore $\omega(\Hf_Y,\cdot)=\mathrm dH_Y$ vanishes on $T_p\D^o$. The conclusion follows by observing that, since $\omega$ is nondegenerate, the $\omega$-orthogonal complement of $T_p\D^o$ in $T_pT\M^*$ has dimension 
$n-k=\mathrm{dim}(\D_x)$.
\end{proof}
\begin{cor}\label{thm:carachamilton}
The projection of a characteristic curve of $\D$ is automatically horizontal. Moreover, let $\gamma:[a,b]\mapsto\M$ be a horizontal curve, let $X_1,\ldots,X_n$ be a time-dependent referential of $T\M$ adapted to $\D$, defined in an open subset $U\subset\R\times\M$ containing the graph of $\gamma$. Define a time-dependent vector field $X$ in $U$ as in \eqref{eq:defX}. Let $\eta:[a,b]\mapsto\D^o$ be a curve with $\pi\circ\eta=\gamma$. Then $\eta$ is a characteristic curve of $\D$ if and only if $\eta$ is an integral curve of $\Hf_X$.
\end{cor}
\begin{proof}
For $p\in\D^o$, the kernel of the restriction of $\omega$ to $T_p\D^o$ is equal to the intersection of $T_p\D^o$ with the $\omega$-orthogonal complement of $T_p\D^o$ in $T_pT\M^*$. By Lemma~\ref{thm:orthcompl}, it follows that the kernel of $\omega\vert_{\D^o}$ projects by $\mathrm d\pi$ into $\D$, and therefore the projection of a characteristic is always horizontal.

For the second part of the statement, observe that for $t\in[a,b]$, $X(t,\cdot)$ is a horizontal vector field in an open neighborhood of $\gamma(t)$ whose value at $\gamma(t)$ is $\dot\gamma(t)$. Therefore $\dot\eta(t)$ is $\omega$-orthogonal to $T_{\eta(t)}\D^o$ if and only if $\dot\eta(t)=\Hf_X(\eta(t))$.
\end{proof}
\begin{cor}\label{thm:sistemaadjunto}
Let $\gamma:[a,b]\mapsto\M$ be a horizontal curve and let $X_1,\ldots,X_n$ be a time-dependent referential of $T\M$ adapted to $\D$, defined in an open subset $U\subset\R\times\M$ containing the graph of $\gamma$. Let $X$ be defined as in \eqref{eq:defX}. A curve $\eta:[a,b]\mapsto\D^o$ with $\pi\circ\eta=\gamma$ is a characteristic of $\D$ if and only if its representation $\tilde\eta(t)\in{\R^n}^*$ in any coordinate chart of $\M$ satisfies the following first order homogeneous linear differential equation:
\begin{equation}\label{eq:sistemaadjunto}
\frac{\mathrm d}{\mathrm dt}\tilde\eta(t)=-\frac{\partial\tilde X}{\partial q}(t,\gamma(t))^*\tilde\eta(t),
\end{equation}
where $\tilde X$ is the representation in coordinates of $X$.
\end{cor}
\begin{proof}
Simply use Corollary~\ref{thm:carachamilton} and write the Hamilton  equations of
$\Hf_X$ in coordinates.
\end{proof}

Differential equation \eqref{eq:sistemaadjunto} is called the {\em adjoint system\/} of \eqref{eq:fundamentaldif}. It is easily seen that $\tilde\eta$ is a solution of \eqref{eq:sistemaadjunto} if and only if $\tilde\eta(t)\tilde v(t)$ is constant for every solution $\tilde v$ of \eqref{eq:fundamentaldif}. From this observation we get:

\begin{lem}\label{thm:transinv}
Let $\gamma:[a,b]\mapsto\M$ be a horizontal curve and suppose that the referential $X_1,\ldots,X_n$ defining $\Phi_t$ in \eqref{eq:defPhi} is adapted to $\D$. Then a curve $\eta:[a,b]\mapsto\D^o$ with $\pi\circ\eta=\gamma$ is a characteristic for $\D$ if and only if $\eta(t)=(\Phi_t^*)^{-1}(\eta(a))$ for every $t\in[a,b]$.
\end{lem}
\begin{proof}
By Corollary~\ref{thm:sistemaadjunto} and the observation above we get that $\eta$ is a characteristic if and only if $\eta(t)v(t)$ is constant for every $v\in\mathrm{Ker}\,\mathrm d\mathcal B_\gamma$. The conclusion follows.
\end{proof}

We can finally prove the main theorem of the section.
\begin{teo}\label{thm:imageend}
The annihilator of the image of the differential of the restriction of the endpoint mapping to $H^1_{q_0}([a,b],\D,\M)$ is given by:
\begin{equation}\label{eq:17}
\begin{split}
\mathrm{Im}\Big(&\mathrm d(\mathrm{end}\vert_{H^1_{q_0}([a,b],\D,\M)})(\gamma)\Big)^o=\\
&\Big\{\eta(b):\eta\ \text{is a characteristic for $\D$ and}\
\pi\circ\eta=\gamma\Big\}
\end{split}
\end{equation}
\end{teo}
\begin{proof}
By Lemma~\ref{thm:tangspace}, we have:
\begin{equation}\label{eq:imend}
\begin{split}
\mathrm{Im}\Big(&\mathrm d(\mathrm{end}\vert_{H^1_{q_0}([a,b],\D,\M)})(\gamma)\Big)=\\
&\left\{\Phi_b\int_a^b\Phi_s^{-1}w(s)\;\mathrm ds:w\ \text{is a $L^2$ horizontal
vector field along $\gamma$}\right\}.
\end{split}
\end{equation}
By Lemma~\ref{thm:transinv}, if $\eta$ is a characteristic with 
$\pi\circ\eta=\gamma$ then $\eta(b)$ annihilates the right hand side of 
\eqref{eq:imend}. 
Namely:
\begin{equation}\label{eq:Phieta}\begin{split}
\eta(b)\Big(\Phi_b\int_a^b\Phi_s^{-1} w(s)& \,\mathrm ds\Big)
=(\Phi_b^*)^{-1}(\eta(a))\left(\Phi_b\int_a^b\Phi_s^{-1} w(s) \,\mathrm
ds\right)\\&=\eta(a)\left(\int_a^b\Phi_s^{-1} w(s) \,\mathrm ds\right)= 
\int_a^b\eta(a)\Phi_s^{-1} w(s) \,\mathrm ds\\&=
\int_a^b(\Phi^*_s)^{-1}\eta(a) w(s)\;\mathrm ds=\int_a^b\eta(s)w(s)\;\mathrm ds
=0.
\end{split} 
\end{equation}
We have to prove that if $\eta_0\in T_{\gamma(b)}\M^*$ annihilates the righthand
side of \eqref{eq:imend} then there exists a characteristic $\eta$ with
$\pi\circ\eta=\gamma$ and $\eta(b)=\eta_0$.

Define $\eta$ by $\eta(t)=(\Phi_t^*)^{-1}(\Phi_b^*(\eta_0))$ for all $t\in[a,b]$. By Lemma~\ref{thm:transinv}, 
we only have to prove that $\eta([a,b])\subset\D^o$. Computing as in
\eqref{eq:Phieta}, we see that, since
$\eta_0$ annihilates the righthand side of \eqref{eq:imend}, then:
\[\int_a^b\eta(s)w(s)\mathrm ds=0,\]
for any horizontal $L^2$ vector field $w$ along $\gamma$. The conclusion follows.
\end{proof}
\begin{cor}\label{thm:imend}
The image of the differential of the restriction of the endpoint mapping 
to $H^1_{q_0}([a,b],\D,\M)$ contains 
$\D_{\gamma(b)}$. 
\end{cor}
\begin{proof}
By Theorem~\ref{thm:imageend}, the annihilator of the image of
the differential of the restriction of the endpoint mapping 
to $H^1_{q_0}([a,b],\D,\M)$ is contained in the annihilator
of $\D_{\gamma(b)}$. The conclusion follows.
\end{proof}
The next corollary, which is obtained easily from
\eqref{eq:17}, gives a characterization of singular curves in terms 
of characteristics: 
\begin{cor}
An $H^1$ curve $\gamma:[a,b]\mapsto\M$ is singular if and only if it is the projection of a non zero characteristic of $\D$.\qed
\end{cor}

Observe that by Lemma~\ref{thm:transinv} a characteristic either never vanishes or is identically zero.
\end{section}

\begin{section}{The normal geodesics as critical points\\ of
the action functional}\label{sec:critical}
In this section we prove that the normal geodesics in $(\M,\D,g)$
correspond to the critical points of the sub-Riemannian action functional
defined in the space of horizontal curves joining two subsets of $\M$.
To this aim, we need to introduce a Lagrangian formalism that will be
be related to the Hamiltonian setup described in Section~\ref{sec:normalgeo}
via the Legendre transform.

We consider the sub-Riemannian action 
functional $E_{\mathrm{sR}}$ defined in the space $H^1([a,b],\D,\M)$:
\begin{equation}\label{eq:defEsR}
E_{\mathrm{sR}}(\gamma)=\frac12\int_a^bg(\dot\gamma,\dot\gamma)\;\mathrm dt.
\end{equation}
The problem of minimizing the action functional $E_{\mathrm{sR}}$ is
essentially equivalent to the problem of minimizing length (see 
Lemma~\ref{thm:energylength} and Corollary~\ref{thm:affsubRiem}).

By Corollary~\ref{thm:Thetasubendsub}, given
$q_0,q_1\in\M$, the set $H^1_{q_0,q_1}([a,b],\D,\M)$
has the structure of a smooth manifold around the regular curves.
It is easy to prove that $E_{\mathrm{sR}}$ is smooth in any open
subset of $H^1_{q_0,q_1}([a,b],\D,\M)$ which has the structure of a
smooth manifold; such an open set will be called a {\em regular\/}
subset of $H^1_{q_0,q_1}([a,b],\D,\M)$. 
We will say that a curve $\gamma\in H^1_{q_0,q_1}([a,b],\D,\M)$
is a critical point of $E_{\mathrm{sR}}$ if  it lies in a regular
subset of $H^1_{q_0,q_1}([a,b],\D,\M)$ and if it is a critical point
of the restriction of $E_{\mathrm{sR}}$ to this regular subset.
The purpose of this section is to prove that the normal geodesics
coincide with the critical points of the $E_{\mathrm{sR}}$ in
$H^1_{q_0,q_1}([a,b],\D,\M)$.

To this goal, we will consider an extension $E$ of $E_{\mathrm{sR}}$ to
the smooth manifold $H^1([a,b],\M)$ defined in terms of
the Riemannian extension $\overline g$ of the sub-Rieman\-nian
metric $g$ that was introduced in Section~\ref{sec:diff}:
\[\phantom{\qquad \gamma\in H^1([a,b],\M)}
E(\gamma)=\frac12\int_a^b\overline g(\dot\gamma,\dot\gamma)\;\mathrm dt,
\qquad \gamma\in H^1([a,b],\M).\]
Let $\gamma\in H^1_{q_0,q_1}([a,b],\D,\M)$ be a regular curve and let
$\theta$ be the map defined in a neighborhood of the graph of $\gamma$
given in Proposition~\ref{thm:kform}. By the method
of Lagrange multipliers, we know that $\gamma$ is a critical
point of $E_{\mathrm{sR}}$ if and only if there exists
$\lambda\in L^2([a,b],\R^n)$ such that $\gamma$ is a critical
point in $H^1_{q_0,q_1}([a,b],\M)$ of the action functional:
\begin{equation}\label{eq:defElambda}
E_\lambda(\gamma)=E(\gamma)-\int_a^b\lambda(t)\cdot\theta_{(t,\gamma(t))}
(\dot\gamma(t))\;\mathrm dt.
\end{equation}
We will see in the proof of Proposition~\ref{thm:normgeocrit} below that
the Lagrange multiplier $\lambda$ associated to a critical point
of $E_{\mathrm{sR}}$ is indeed a smooth map.

$E_\lambda$ is the action functional of the time-dependent
Lagrangian $\mathcal L_\lambda$ defined on an open subset of $T\M$, given by:
\begin{equation}\label{eq:defLlalmbda}
\mathcal L_\lambda(t,v)=\frac12 \overline g(v,v)-\lambda(t)\cdot
\theta_{(t,m)}(v),\quad v\in T_m\M.
\end{equation}
The Lagrangian $\mathcal L_\lambda$ is $L^1$ in the variable $t$,
moreover,  for (almost) all $t\in[a,b]$, the map $v\mapsto \mathcal L_{\lambda}(t,v)$
is smooth. Therefore the critical points of $E_\lambda$ are 
curves satisfying the Euler--Lagrange equations; in a chart
$q=(q_1,\ldots,q_n)$, the equations are:
\begin{equation}\label{eq:EL}
\frac{\partial \mathcal L_\lambda}{\partial q}-\frac{\mathrm d}{\mathrm dt}
\,\frac{\partial \mathcal L_\lambda}{\partial\dot q}=0.
\end{equation}
We recall that if $\mathcal L:U\subset\R\times T\M$ is a time-dependent
Lagrangian defined on an open subset of $\R\times T\M$, the 
{\em fiber derivative\/} of $\mathcal L$ is the map
$\mathbb F\mathcal L:U\mapsto\R\times T\M^*$ given by:
\[\mathbb F\mathcal L(t,v)=(t,\mathrm d(\mathcal L\vert_{U\cap T_{\pi(v)}\M})
(v)),\]
where $\pi:T\M\mapsto\M$ is the projection. 
For $t\in\R$, we denote by $U_t$ the open subset of $T\M$
consisting of those $v$'s such that $(t,v)\in U$.
The Lagrangian $\mathcal L$ is said to be {\em regular\/} if, 
for each $t$,
the map  $v\mapsto\mathbb F\mathcal L(t,v)$
is a local diffeomorphism; $\mathcal L$ is said to be {\em hyper-regular\/}
if $v\mapsto\mathbb F\mathcal L(t,v)$ is a diffeomorphism between $U_t$ 
and an open subset of $T\M^*$. Associated to a hyper-regular Lagrangian
$\mathcal L$ in $U\subset\R\times T\M$ 
one has a Hamiltonian $H$ defined on the open subset
$\mathbb F\mathcal L(U)$ by the formula:
\[\phantom{\quad (t,v)\in U.}
H\left(\mathbb F\mathcal L(t,v)\right)=\mathbb FL(t,v)v-\mathcal L(t,v),
\quad (t,v)\in U.\]
This procedure is called the {\em Legendre transform\/}
(see \cite[Chapter~3]{AM}).
If $\mathcal L$ is a hyper-regular Lagrangian and
$H$ is the associated Hamiltonian, then the solutions of the 
Euler--Lagrange equations \eqref{eq:EL}
of $\mathcal L$ correspond, via $\mathbb F\mathcal L$,
to the solutions of the Hamilton  equations of $H$, i.e.,
a smooth curve $\gamma:[a,b]\mapsto\M$ is a solution of
\eqref{eq:EL} if and only if $\Gamma=\mathbb F\mathcal L\circ(\gamma,\dot\gamma)$
is a solution of the Hamiltonian $H$.

Let us show now the this formalism applies to the case of
the Lagrangian $\mathcal L_\lambda$ of \eqref{eq:defLlalmbda}:
\begin{lem}\label{thm:Llambdahyperreg}
The Lagrangian $\mathcal L_\lambda$ is hyper-regular.
\end{lem}
\begin{proof}
From \eqref{eq:defLlalmbda},
the fiber derivative $\mathbb F\mathcal L_\lambda$ is easily computed
as:
\begin{equation}\label{eq:derLlambda}
\mathbb F\mathcal L_\lambda(t,v)=\overline g(v,\cdot)-\lambda(t)\cdot\theta_{(t,m)}
\in T_m\M^*.
\end{equation}
For each $t\in[a,b]$, the map $\mathbb F\mathcal L_\lambda(t,\cdot):
T_m\M\mapsto T_m\M^*$ is clearly a diffeomorphism, whose inverse
is given by:
\begin{equation}\label{eq:using}
T_m\M^*\ni p\mapsto \overline g^{-1}\left(p+\lambda(t)\cdot\theta_{(t,m)}\right)
\in T_m\M.
\end{equation}
\end{proof}
We are finally ready to prove the following:
\begin{prop}\label{thm:normgeocrit}
Let $\gamma$ be a regular curve in $H^1_{q_0,q_1}([a,b],\D,\M)$.
Then, $\gamma$ is a critical point of $E_{\mathrm{sR}}$ if and only if
it is a normal sub-Riemannian geodesic in $(\M,\D,g)$.
\end{prop}
\begin{proof}
A critical point of $E_{\mathrm{sR}}$ is a  curve
satisfying the Euler--Lagrange equations \eqref{eq:EL}
associated to  the Lagrangian  $\mathcal L_\lambda$ of \eqref{eq:defLlalmbda}.
By Lemma~\ref{thm:Llambdahyperreg}, $\mathcal L_\lambda$ is hyper-regular,
hence the solutions of \eqref{eq:EL} correspond, via $\mathbb F\mathcal L_\lambda$
to the solutions of the associated  Hamiltonian $H_\lambda$, computed as follows.
First, for $v\in T_m\M$ we have:
\[\begin{split}
\mathbb F\mathcal L_\lambda&(t,v)\,v-\mathcal L_\lambda(t,v)=\\ &=
\overline g(v,v)-\lambda(t)\cdot\theta_{(t,m)}(v)-\frac12\,\overline g(v,v)+\lambda(t)\cdot\theta_{(t,m)}(v)
=\frac12\,\overline g(v,v).
\end{split}
\]
Then, using \eqref{eq:using}, we compute:
\begin{equation}\label{eq:Hlambda}
H_\lambda(t,q,p)=\frac12\,\overline g^{-1}\big(p+\lambda(t)\cdot\theta_{(t,q)},\,
p+\lambda(t)\cdot\theta_{(t,q)}\big).
\end{equation}
For the proof of the Proposition, we need to show that if $\gamma$ is 
an absolutely continuous  curve in $\M$, then $\gamma$ is horizontal and it is a solution for the 
Hamilton 
equations associated to the Hamiltonian $H_\lambda$ for some
$\lambda$ if and only if
it is a solution of the Hamilton  equations associated to the
sub-Riemannian Hamiltonian $H$ of formula \eqref{eq:defHsR}.

The Hamilton  equations of $H_\lambda$ are computed as follows:
\begin{equation}\label{eq:HJHlambda}
\left\{\begin{array}{l} \displaystyle
\frac{\mathrm dq}{\mathrm dt}=\overline g^{-1}(p+\lambda(t)\cdot\theta_{(t,q)});
\\ \\ \displaystyle \frac{\mathrm dp}{\mathrm dt}=-\overline g^{-1}\big(
\lambda(t)\cdot\frac{\partial\theta_{(t,q)}}{\partial q},p+\lambda(t)
\cdot\theta_{(t,q)}\big).
\end{array}\right.
\end{equation} 
From the horizontality of $\frac{\mathrm dq}{\mathrm dt}$, using the first
equation of \eqref{eq:HJHlambda} we get:
\[\Big(p+\lambda(t)\cdot\theta_{(t,q)}\Big)\big\vert_{\D_1}=0,\]
and since $\theta\vert_{\D_1}$ is an isomorphism, we get an explicit
expression for the Lagrange multiplier $\lambda$:
\begin{equation}\label{eq:explicitlambda}
\lambda(t)=-p(t)\circ\Big[\theta_{(t,q)}\big\vert_{\D_1}\Big]^{-1}.
\end{equation}
Observe that, by a standard {\em boot-strap\/} argument, from
\eqref{eq:explicitlambda} it follows easily that $\lambda$
is smooth.

We now write the Hamilton equations of the sub-Riemannian Hamiltonian
and of $H_\lambda$ using a suitable time-dependent
referential $X_1,\ldots,X_n$ of $T\M$. The choice of the referential is done as follows.
Let $\theta_1,\ldots,\theta_k$ be a time-dependent referential of the annihilator
$\D^o=\left(\D^\perp\right)^*$ which is orthonormal with respect to
$\overline g^{-1}$. For the orthogonality, it suffices to consider any referential
of $\D^o$ and then to orthonormalize it by the method of Gram-Schmidt.
Then, let $X_{n-k+1},\ldots,X_n$ be the  referential of $\D^\perp$ obtained
by dualizing $\theta_1,\ldots,\theta_k$. Finally, let $X_1,\ldots,X_{n-k}$
be any orthonormal referential of $\D$, time-dependent or not.

In the referential $X_1,\ldots,X_n$, for $i=1,\ldots,n-k$ we have:
\begin{equation}\label{eq:intheref}
\left[\theta_{(t,q)}\big\vert_{\D_1}\right]^{-1}\left[\frac{\partial\theta_{(t,q)}}{%
\partial q}\,X_i\right]=\sum_{j=1}^k\left[\frac{\partial \theta_j}{\partial q}
(t,q)\,X_i\right]\cdot X_{n-k+j}.
\end{equation}
We can rewrite \eqref{eq:HJHlambda} as:
\begin{equation}\label{eq:newHJHlambda}
\left\{\begin{array}{l}\displaystyle \frac{\mathrm dq}{\mathrm dt}=
\sum_{i=1}^{n-k}p(X_i)\,X_i+\sum_{i=n-k+1}^n\left(p(X_i)+\lambda_{i-n+k}\right)\,X_i,\\ \\
\displaystyle \frac{\mathrm dp}{\mathrm dt}=-\sum_{i=1}^{n-k}p(X_i)\,p\left(
\frac{\partial X_i}{\partial q}\right)-\sum_{i=n-k+1}^n2\left(p(X_i)+\lambda_{i-n+k}\right)
\,p\left(\frac{\partial X_i}{\partial q}\right),
\end{array}\right.
\end{equation}
where $\lambda=(\lambda_1,\ldots,\lambda_k)$. On the other hand, the Hamilton 
equations for $H$ are written as:
\begin{equation}\label{eq:newHJH}
\left\{\begin{array}{l}\displaystyle \frac{\mathrm dq}{\mathrm dt}=
\sum_{i=1}^{n-k}p(X_i)\,X_i,\\ \\
\displaystyle \frac{\mathrm dp}{\mathrm dt}=-\sum_{i=1}^{n-k}p(X_i)\;p\!\left(
\frac{\partial X_i}{\partial q}\right).
\end{array}\right.
\end{equation}
Now, if $\gamma$ is horizontal and it satisfies \eqref{eq:newHJHlambda} for some
$\lambda$ it follows that the second sum of the first equation in \eqref{eq:newHJHlambda}
is zero, and therefore $\gamma$ satisfies also \eqref{eq:newHJH}.
Conversely, if $\gamma$ satisfies \eqref{eq:newHJH}, then $\gamma$ is horizontal,
and defining $\lambda$ by \eqref{eq:explicitlambda}, it is easily seen that
$\gamma$ is a solution of \eqref{eq:HJHlambda}.
\end{proof}
\begin{rem}\label{thm:remlift}
It follows easily from \eqref{eq:newHJH} that if $\gamma$ is a normal
geodesic and $\Gamma$ is a Hamiltonian lift of $\gamma$, then
$\Gamma\vert_{\mathcal D}=g(\dot\gamma,\cdot)$.
\end{rem}
We now consider the case of sub-Riemannian geodesics with endpoints varying 
in two submanifolds of $\M$.
\begin{prop}\label{thm:critpts}
Let $(\M,\D,g)$ be a sub-Riemannian manifold, let $\mathcal P,\mathcal Q\subset\M$ be
smooth submanifolds of $\M$ and assume that $\mathcal Q$ is {\em transversal\/}
to $\D$, i.e., $T_q\mathcal Q+\D_q=T_q\M$ for all $q\in \mathcal Q$.
Then, the set
\[H^1_{\mathcal P,\mathcal Q}([a,b],\D,\M)=\Big\{x\in H^1([a,b],\D,\M):
x(a)\in\mathcal P,\ x(b)\in\mathcal Q\Big\}\]
is a smooth submanifold of $H^1([a,b],\M)$. Moreover, the critical points
of the sub-Riemannian action functional $E_{\mathrm{sR}}$ in 
$H^1_{\mathcal P,\mathcal Q}([a,b],\D,\M)$
are precisely the normal geodesics $\gamma$ joining 
$\mathcal P$ and $\mathcal Q$
that admit a lift $\Gamma:[a,b]\mapsto T\M^*$ satisfying 
the boundary conditions:
\begin{equation}\label{eq:lift}
\Gamma(a)\in T_{\gamma(a)}\mathcal P^o,\quad\text{and}\quad
\Gamma(b)\in T_{\gamma(b)}\mathcal Q^o.
\end{equation}
\end{prop}
\begin{proof}
The fact that $H^1_{\mathcal P,\mathcal Q}([a,b],\D,\M)$
is a smooth manifold follows easily from the transversality
of $\mathcal Q$ and Corollary~\ref{thm:imend}.

The proof of the second part of the statement is analogous to
the proof of Proposition~\ref{thm:normgeocrit}, keeping in mind
that the critical points of the action functional associated to
a hyper-regular Lagrangian in the space of  curves joining
$\mathcal P$ and $\mathcal Q$ are the solutions of the Hamilton 
equations whose Hamiltonian lift vanishes on the tangent spaces
of $\mathcal P$ and $\mathcal Q$.
\end{proof}
Obviously, the role of $\mathcal P$ and $\mathcal Q$ in
Proposition~\ref{thm:critpts} can be interchanged, and the
same conclusion holds in the case that $\mathcal P$ is 
transversal to $\D$.

As a consequence of Proposition~\ref{thm:critpts} we get some information
on the geometry of sub-Riemannian {\em balls}. Given a horizontal
curve $\gamma:[a,b]\mapsto\M$, we define $\ell(\gamma)$ to be its length:
\[\ell(\gamma)=\int_a^bg(\dot\gamma,\dot\gamma)^{\frac12}\,\mathrm dt.\]  
For $q_0,q_1\in \M$, we set
\[\mathrm{dist}(q_0,q_1)=\inf\Big\{\ell(\gamma):\gamma\ \text{is a
horizontal curve joining} \ q_0\ \text{and}\ q_1\Big\}\in[0,+\infty],\]
where such number is infinite if $q_0$ and $q_1$ cannot be joined by
any horizontal curve. 
A horizontal curve $\gamma:[a,b]\mapsto\M$ is said to be {\em length
minimizing\/} between two subsets $\mathcal P$ and $\mathcal Q$ of $\M$ if
$\gamma(a)\in\mathcal P$, $\gamma(b)\in\mathcal Q$ and
\[\ell(\gamma)=\inf_{\stackrel{\scriptstyle q_0\in\mathcal P}{%
q_1\in\mathcal Q}}\mathrm{dist}(q_0,q_1).\] A horizontal
curve
$\gamma$ is said to be {\em affinely parameterized\/} if
$g(\dot\gamma,\dot\gamma)$ is almost everywhere constant. Every horizontal
curve is the reparameterization of an affinely parameterized horizontal curve
(see Corollary~\ref{thm:affsubRiem}). Since the sub-Riemannian Hamiltonian
is constant on its integral curves, it follows that every
normal geodesic is affinely parameterized. Moreover, using the
Hamilton  equations \eqref{eq:newHJH}, it is easy to see that
an affine reparameterization of a normal geodesic is again a normal
geodesic.

We relate the problem of minimization of the length and of the action
functional by the following:
\begin{lem}\label{thm:energylength}
Let $\gamma:[a,b]\mapsto\M$ be an horizontal curve joining
the submanifolds $\mathcal P$ and $\mathcal Q$. Then, $\gamma$ is a minimum
of $E_{\mathrm{sR}}$ in $H^1_{\mathcal P,\mathcal Q}([a,b],\D,\M)$
if and only if $\gamma$ is affinely parameterized and $\gamma$ is a length
minimizer between $\mathcal P$ and $\mathcal Q$.
\end{lem}
\begin{proof}
By Cauchy--Schwartz inequality we have:
\[\ell(\gamma)^2\le2(b-a) E_{\mathrm{sR}}(\gamma)^2, \]
where the equality holds if and only if $\gamma$ is affinely parameterized.
If $\gamma$ is affinely parameterized and it minimizes length, then, for any
$\mu\in H^1_{\mathcal P,\mathcal Q}([a,b],\D,\M)$, we have:
\[E_{\mathrm{sR}}(\gamma)=\frac{\ell(\gamma)^2}{2(b-a)}\le 
\frac{\ell(\mu)^2}{2(b-a)}\le E_{\mathrm{sR}}(\mu).\]
Hence, $\gamma$ is a minimum of $E_{\mathrm{sR}}$.

Conversely, suppose that $\gamma$ is a minimum of $E_{\mathrm{sR}}$.
There exists an affinely parameterized horizontal curve $\mu:[a,b]\mapsto
\M$ such that $\gamma$ is a reparameterization of $\mu$ (see
Corollary~\ref{thm:affsubRiem}). We have:
\[E_{\mathrm{sR}}(\gamma)\le E_{\mathrm{sR}}(\mu)=
\frac{\ell(\mu)^2}{2\,(b-a)}=\frac{\ell(\gamma)^2}{2\,(b-a)}\le
E_{\mathrm{sR}}(\gamma),\]
hence the above inequalities are indeed equalities, and $\gamma$ is
affinely parameterized. 

Now, assume by contradiction that $\rho:[a,b]\mapsto\M$ connects 
$\mathcal P$ and $\mathcal Q$ and satisfies $\ell(\rho)<\ell(\gamma)$.
By Corollary~\ref{thm:affsubRiem}, we can assume that $\rho$ is affinely
parameterized, hence $E_{\mathrm{sR}}(\rho)< E_{\mathrm{sR}}(\gamma)$.
This is a contradiction, and we are done.
\end{proof}

For $q_0\in\M$ and $r\in\R^+$, the {\em open ball\/} $B_r(q_0)$
is defined by:
\[B_r(q_0)=\Big\{q_1:\mathrm{dist}(q_0,q_1)< r\Big\}.\]
\begin{cor}\label{thm:corballs}
Suppose that there exists an affinely parameterized length minimizer 
$\gamma:[a,b]\mapsto\M$ between
$q_0$ and $q_1$  which is not a normal extremal; set $r=\mathrm{dist}(q_0,q_1)$.
Then, any submanifold $\mathcal Q$ through $q_1$  which is transversal to
$\D$ at $q_1$ has non empty intersection with the open ball
$B_r(q_0)$.
\end{cor}
\begin{proof}
By contradiction, suppose that we can find a submanifold $\mathcal Q$
through $q_1$ which is transversal to $\D$ at $q_1$ and disjoint from
the open ball $B_r(q_0)$. It follows that $\gamma$ is a length minimizer
between the point $q_0$ and the submanifold $\mathcal Q$, hence, by 
Lemma~\ref{thm:energylength},
$\gamma$ is a minimum point for the action functional in 
$H^1_{q_0,\mathcal Q}([a,b],\D,\M)$. By possibly considering a small
portion of $\mathcal Q$ around $q_1$, we can assume that $\mathcal Q$
is everywhere transversal to $\D$. From Proposition~\ref{thm:critpts} it follows
then that $\gamma$ is a normal geodesic, which is a contradiction.
\end{proof}
\begin{rem}\label{thm:remsmoothness}
Proposition~\ref{thm:critpts} can also be used to establish
the {\em smoothness\/} of abnormal extremizers, which is in general an open
question. 
Observe indeed that its statement can be rephrased as follows.
Let $\gamma\colon[a,b]\to\M$ be an affinely parameterized
length-minimizer connecting $q_0$ and $q_1$ in $\M$; set
$r=\mathrm{dist}(q_0,q_1)$. If there exists a manifold $\mathcal Q$ transverse to
$\mathcal D$ passing through $q_1$ which does not intercept the open ball
$B(q_0;r)$ then $\gamma$ is a normal extremal and consequently it is
smooth.
\end{rem}

As a corollary of Proposition~\ref{thm:normgeocrit}, we also obtain
an alternative proof of a result of \cite{LS} that gives necessary conditions
for length minimizing:
\begin{cor}\label{thm:normalabnormal}
An affinely parameterized length minimizer is either an abnormal
minimizer or a normal geodesic.
\end{cor}
\begin{proof}
It follows immediately from Definition~\ref{thm:defregular},
Proposition~\ref{thm:normgeocrit} and the fact that affinely parameterized
length minimizers are minima of the sub-Riemannian action functional.
\end{proof}

The solutions of sub-Riemannian geodesic problem with variable endpoints
in the case that the end-manifold is one-dimensional has a physical
interpretation in the context of General Relativity (see \cite{GP,GPV}). 
Such geodesics  can be interpreted as the solution of a general
relativistic {\em brachistochrone problem\/} in a stationary
Lorentzian manifold.
\end{section}

%%%%%%%%%%%%%%%
%%%%%%%%%%%%%%%
\appendix

\begin{section}{Affine parameterization of horizontal curves}
\label{sec:affine}
In this  appendix we show that every horizontal curve
in a sub-Riemannian manifold can be obtained as the reparameterization
of an affinely parameterized horizontal curve.

Given two absolutely continuous curves $\gamma:[a,b]\mapsto\M$ and
$\mu:[c,d]\mapsto\M$, we say that $\gamma$ is a {\em reparameterization\/}
of $\mu$ if there exists an absolutely continuous, nondecreasing and surjective
map $\sigma:[a,b]\mapsto[c,d]$ such that $\gamma=\mu\circ\sigma$.
It can be proven that in this case
$\dot\gamma=(\dot\mu\circ\sigma)\,\dot\sigma$ almost everywhere.
\begin{prop}\label{thm:affRiem}
Let $(\M,\overline g)$ be a Riemannian manifold,  $\gamma:[a,b]\mapsto\M$ an
absolutely continuous curve. Then, there exists a unique pair
of absolutely continuous maps $\mu:[0,L]\mapsto\M$ and $\sigma:[a,b]\mapsto
[0,L]$, with $\sigma$ nondecreasing and surjective, such that
$\overline g(\dot\mu(t),\dot\mu(t))\equiv1$ almost everywhere on $[0,L]$
and $\gamma=\mu\circ\sigma$.
\end{prop}
\begin{proof}
Suppose that the pair $\mu,\sigma$ satisfying the thesis
is found; then we obtain easily 
\begin{equation}\label{eq:defsigma}
\sigma(t)=\ell(\gamma\vert_{[a,t]})=\int_a^t\overline g(\dot\gamma,\dot\gamma)
^\frac12\;\mathrm dt. 
\end{equation}
Since $\sigma$ is surjective, this proves the uniqueness of the pair.

As to the existence, set $L=\ell(\gamma)$ and define $\sigma$ as in
\eqref{eq:defsigma}. Obviously, $\sigma$ is absolutely continuous,
nondecreasing and surjective. 

Suppose that $\sigma(s)=\sigma(t)$ for some $s,t\in[a,b]$, with $s<t$.
Then, $\ell(\gamma\vert_{[s,t]})=0$, and therefore $\gamma(s)=\gamma(t)$.
It follows that there exists a function $\mu:[0,L]\mapsto\M$ with
$\mu\circ\sigma=\gamma$. The curve $\mu$ is Lipschitz continuous, hence
absolutely continuous; for, if $s,t\in [0,L]$, let
$s_1,t_1\in[a,b]$ be such that $\sigma(s_1)=s$ and $\sigma(t_1)=t$.
Then,
\[\mathrm{dist}(\mu(s),\mu(t))=\mathrm{dist}(\gamma(s_1),\gamma(t_1))\le
\ell(\gamma\vert_{[s_1,t_1]})=\vert\sigma(s_1)-\sigma(t_1)\vert=\vert
s-t\vert.\]
We are left with the proof that $\overline g(\dot\mu,\dot\mu)\equiv1$
almost everywhere. To see this, let $t\in[0,L]$ be chosen and let
$t_1\in[a,b]$ be such that $t=\sigma(t_1)$. Then, we have:
\begin{equation}\label{eq:aserderivada}
\int_0^t\overline
g(\dot\mu,\dot\mu)^\frac12\;\mathrm
dr=\ell(\mu\vert_{[0,t]})=\ell(\gamma\vert_{[a,t_1]})=\sigma(t_1)=t.
\end{equation}
The conclusion follows by differentiating \eqref{eq:aserderivada} 
with respect to $t$.
\end{proof}

\begin{lem}\label{thm:muhorizontal}
Let $\M$ be a smooth manifold and $\D\subset T\M$ be a smooth distribution.
Let $\mu:[a,b]\mapsto\M$ be an absolutely continuous curve; if $\mu$ admits
a reparameterization which is horizontal, then $\mu$ is horizontal.
\end{lem}
\begin{proof}
Let $\sigma:[c,d]\mapsto[a,b]$ an absolutely continuous nondecreasing surjective
map with $\gamma=\mu\circ\sigma$ horizontal. Define:
\[\begin{split}
X=&\;\Big\{t\in[c,d]:\text{the
equality}\ \dot\gamma(t)=\dot\mu(\sigma(t))\dot\sigma(t)\
\text{fails to hold}\Big\},
\\ Y=&\;\Big\{t\in[c,d]:\dot\sigma(t)=0\Big\}.
\end{split}\]
Clearly, $\mu$ is horizontal outside $\sigma(X\cup Y)$; to conclude the
proof it suffices to show that $\sigma(X\cup Y)$ has null measure.
To see this, observe that $X$ has null measure and therefore
$\sigma(X)$ has null measure. Moreover, since $\dot\sigma=0$ in $Y$,
it is not difficult to show that $\sigma(Y)$ has null measure, and we are done. 
\end{proof}

\begin{cor}\label{thm:affsubRiem}
Let $(\M,\D,g)$ be a sub-Riemannian manifold and $\gamma$ a horizontal
curve in $\M$. Then, $\gamma$ is the reparameterization of a unique
horizontal curve $\mu:[0,L]\mapsto\M$ such that $g(\dot\mu,\dot\mu)\equiv1$
almost everywhere.
\end{cor}
\begin{proof}
Let $\overline g$ be any Riemannian extension of $g$ and apply
Proposition~\ref{thm:affRiem}. The curve $\mu$ thus obtained is horizontal
by Lemma~\ref{thm:muhorizontal}.
\end{proof}
\end{section}
%%%%%%%%%%%%%%%%%%%
\begin{section}{Local minimality of normal geodesics}
\label{sec:localminimality}
The aim of this section is to prove that a sufficiently small
segment of a sub-Riemannian normal geodesic is a distance minimizer 
between an initial submanifold and a point. We will simply adapt the 
proof of local optimality presented in \cite[Appendix~C]{LS}.

\begin{prop}\label{thm:locmin}
Let $(\M,\mathcal D,g)$ be a sub-Riemannian manifold, $\mathcal P\subset
\M$ a submanifold and $\gamma:[a,b]\mapsto\M$ a normal geodesic with
$\gamma(a)\in\mathcal P$ and such that there exists
a Hamiltonian lift  $\Gamma:[a,b]\mapsto T\M^*$ of $\gamma$ with
$\Gamma(a)\vert_{T_{\gamma(a)}\mathcal P}=0$. Then, for $\varepsilon>0$ small
enough,
$\gamma\vert_{[a,a+\varepsilon]}$ is a length minimizer between $\mathcal P$
and $\gamma(a+\varepsilon)$.
\end{prop}
\begin{proof}
We can assume without loss of generality that $g(\dot\gamma,\dot\gamma)=1$.
Let $\mathcal S\subset\M$ be a codimension $1$ submanifold containing
a neighborhood of $\gamma(a)$ in $\mathcal P$ and such that
$\Gamma(a)\vert_{T_{\gamma(a)}\mathcal S}=0$. The existence of such
a submanifold is easily proved using a coordinate system in $\M$
adapted to $\mathcal P$ around $\gamma(a)$. Observe that, by
Remark~\ref{thm:remlift}, we have $g^{-1}\big(\Gamma(a)\vert_{\mathcal
D},\Gamma(a)\vert_{\mathcal D}\big)=1$.

Let $\lambda:\mathcal S\mapsto T\M^*$ be a $1$-form in $\M$ along $\mathcal S$
such that $\lambda(x)\vert_{T_x\mathcal S}=0$,
$g^{-1}\big(\lambda(x)\vert_{\mathcal D}, \lambda(x)\vert_{\mathcal D}\big)=1$
for all
$x\in\mathcal S$ and such that $\lambda(\gamma(a))=\Gamma(a)$. Let
$U\subset\mathcal S$ be a sufficiently small open subset
containing $\gamma(a)$ and let $\varepsilon>0$ be sufficiently small. 
Consider the map
$\Phi:\, ]\,a-\varepsilon,a+\varepsilon\,[ \times U\mapsto T\M^*$ 
such that $t\mapsto\Phi(t,x)$ is a solution of the sub-Riemannian
Hamiltonian $H$ defined in \eqref{eq:defHsR} and $\Phi(a,x)=\lambda(x)$ for all
$x\in U$. Let $F=\pi\circ\Phi$, where $\pi:T\M^*\mapsto\M$ is
the projection. 

\noindent
By Remark~\ref{thm:remlift}, $\Gamma(a)(\dot\gamma(a))=1$, which
implies that $T_{\gamma(a)}\M=T_{\gamma(a)}\mathcal
S\oplus(\R \dot\gamma(a))$. It follows easily that the
differential of $F$ at $(a,\gamma(a))$ is an isomorphism, and
by the Inverse Function Theorem, by possibly passing to smaller
$\varepsilon$ and $U$,  $F$ is a diffeomorphism between
$]\,a-\varepsilon,a+\varepsilon\,[\times U$ and an open neighborhood
$V$ of $\gamma(a)$ in $\M$. By possibly taking a smaller $V$, we can 
assume that $V\cap\mathcal P\subset\mathcal S$.

We define a vector field $X$, a $1$-form $\lambda$ and a smooth
map $\tau$ on $V$ by setting:
\[\tau\big(F(t,x)\big)=t,\quad X\big(F(t,x)\big)=\frac{\mathrm d}{\mathrm
dt}\,F(t,x),
\quad\lambda\big(F(t,x)\big)=\Phi(t,x),\]
for all $(t,x)\in\,]\,a-\varepsilon,a+\varepsilon\,[\times U$.
Since $H\circ\Phi$ does not depend on $t$, it follows easily
that 
\begin{equation}\label{eq:normalambda}
g^{-1}\big(\lambda\vert_{\mathcal
D},\lambda\vert_{\mathcal D}\big)=1.
\end{equation}

We prove next that $\lambda=\mathrm d\tau$. To this aim, let $\Psi_X$
denote the flow of $X$, defined on an open subset of $\R\times V$;
for $s\in\R$ we set $\Psi_X^s=\Psi_X(s,\cdot)$. Clearly,
$t\mapsto F(t,x)$ is an integral curve of $X$, and therefore we
have $\tau\circ\Psi_X^s=s+\tau$, hence $\mathrm d\tau$ is invariant
by the flow of $X$, i.e., 
\[(\Psi_X^s)^*(\mathrm d\tau)=\mathrm d\tau.\]
We show that $\lambda$ is also invariant by the flow of $X$;
the equality $\lambda=\mathrm d\tau$ will follow from the fact that
these two $1$-forms coincide on $\mathcal S$.
For the invariance of $\lambda$, we argue as follows: let $x\in U$,
$v_0\in T_x\M$ and $v(t)=\mathrm d\Psi_X^{t-a}(x)[v_0]$; it suffices to prove
that $\lambda(F(t,x))(v(t))$ is constant in $t$.

In local coordinates $q=(q_1,\ldots,q_n)$, $v$ satisfies the
following linear differential equation:
\begin{equation}\label{eq:eqdiffv}
\frac{\mathrm dv}{\mathrm dt}=\frac{\partial X}{\partial q}(v).
\end{equation}
For $t\in\,]-\varepsilon,\varepsilon[$ fixed, let
$X_1,\ldots,X_{n-k}$ be an orthonormal frame for $\mathcal D$ around
$F(t,x)$; by Remark~\ref{thm:remlift} we have
$\Phi(t,x)\vert_{\mathcal D}=g(X(F(t,x)),\cdot)$, from which it
follows:
\begin{equation}\label{eq:Xlambda}
X=\sum_{i=1}^{n-k}\lambda(X_i)\,X_i.
\end{equation}
From \eqref{eq:normalambda} it follows that $\sum_i\lambda(X_i)^2=1$, and
differentiating this expression we obtain:
\begin{equation}\label{eq:sum}
\sum_{i=1}^{n-k}\lambda(X_i)\,\frac{\partial}{\partial
q}\left(\lambda(X_i)\right)=0.
\end{equation}
From \eqref{eq:Xlambda} and \eqref{eq:sum}, it follows:
\begin{equation}\label{eq:lambdadelX}
\lambda\left(\frac{\partial X}{\partial q}\right)=\sum_{i=1}^{n-k}
\lambda(X_i)\,\lambda\left(\frac{\partial X_i}{\partial q}\right).
\end{equation}
Using the second Hamilton  equation in \eqref{eq:newHJH}, we finally get:
\begin{equation}\label{eq:conclus}
\frac{\mathrm d}{\mathrm
dt}\,\lambda\big(F(t,x)\big)=-\lambda\left(\frac{\partial X}{\partial
q}\right).
\end{equation}
Using \eqref{eq:eqdiffv} and \eqref{eq:conclus} it is easily seen
that $\lambda(F(t,x))v(t)$ is constant in $t$, and $\lambda$ is
invariant by the flow of $X$.

The equality $\lambda=\mathrm d\tau$ is thus proven, and by
\eqref{eq:normalambda} we obtain:
\begin{equation}\label{eq:normadtau}
g^{-1}\big(\mathrm d\tau\vert_{\mathcal
D},\mathrm d\tau\vert_{\mathcal D}\big)=1.
\end{equation}
Let now $\mu:[a,a+\varepsilon]\mapsto V$ be a horizontal curve
with $\mu(a)\in\mathcal P$ and $\mu(a+\varepsilon)=\gamma(a+\varepsilon)$.
Using \eqref{eq:normadtau}, the length of $\mu$ is estimated as follows:
\[\begin{split}L(\mu)=\int_a^{a+\varepsilon}\Vert\dot\mu\Vert\,\mathrm dt&\ge
\int_a^{a+\varepsilon}\mathrm d\tau(\dot\mu(t))\;\mathrm
dt=\tau(\mu(a+\varepsilon))-\tau(\mu(a))=\varepsilon=\\ \\ &
=L\big(\gamma\vert_{[a,a+\varepsilon]}\big).
\end{split}\]
This implies that $\gamma\vert_{[a,a+\varepsilon]}$ is a length minimizer
between $\mathcal P$ and $\gamma(a+\varepsilon)$ among all the horizontal
curves with image in $V$. The conclusion of the proof will follow
from the next Lemma, by possibly considering a smaller $\varepsilon$.
\end{proof}
\begin{lem}\label{thm:lengthminloc}
Let $(\M,\D,g)$ be a sub-Riemannian manifold and let $V\subset\M$ 
be an open subset. Given $x\in U$ there exists $r>0$ such that
every horizontal curve $\mu:[a,b]\mapsto\M$ with $\mu(a)=x$ and
$L(\mu)<r$ satisfies $\mu([a,b])\subset V$.
\end{lem}
\begin{proof} 
We compare the sub-Riemannian metric $g$ with the Euclidean metric relative to
an arbitrary coordinate system around $x$. Let $\varphi:W\mapsto\widetilde W$
be a coordinate system in $\M$ with $x\in W$, $W\subset V$ and $\widetilde
W$ is an open neighborhood of $0$ in $\R^n$. Let
$B\subset W$ be the inverse image through $\varphi$ of a closed
ball of radius $s$, $B[\varphi(x);s]\subset\widetilde W$. For $m\in W$
and  $v\in T_m\M$, denote by $\Vert v\Vert_e$ the Euclidean norm
of the vector $\mathrm d\phi(m)[v]$. The set of vectors
$v\in\D$ that are tangent to the points of $B$ with $\Vert v\Vert_e=1$
form a compact subset of $T\M$, in which the continuous function 
$v\mapsto g(v,v)^{\frac12}=\Vert v\Vert$ attains a positive minimum $k$.
Observe that for all $v\in\D$ tangent to some point of $B$, it is
$\Vert v\Vert\ge k\cdot\Vert v\Vert_e$.

Take $r=ks>0$. If $\mu:[a,b]\mapsto\M$ is a horizontal curve with $\mu(a)=x$
and $\mu([a,b])\not\subset V$, then there exists $c\in\,]a,b\,[$ with
$\mu([a,c]\subset B$ and $\gamma(c)\in\partial B$. Therefore, 
\[L(\mu)\ge L(\mu\vert_{[a,c]})\ge k L_e\big(\varphi\circ\mu\vert_{[a,c]}\big)\ge
ks=r,\]
where $L_e$ denotes the Euclidean length of a curve. This concludes the proof.
\end{proof}
\end{section}


\begin{thebibliography}{99}

\bibitem{AM} R.\ Abraham, J.\ E.\ Marsden, {\em Foundations of Mechanics}, 
2nd Edition, Benjamin/Cummings, Advanced Book Program, Reading, Mass., 1978.

\bibitem{Bis} J.-M.\ Bismut, {\em Large Deviations and the Malliavin
Calculus}, Progress in Mathematics, Birkh\"auser, Boston - Basel - Stuttgart, 1984.

\bibitem{dC} M.\ do Carmo, {\em Riemannian Geometry}, 
Birkh\"auser, Boston, 1992.

\bibitem{GGP} R.\ Giamb\'o, F.\ Giannoni, P.\ Piccione,
{\em Existence, Multiplicity and Regularity for sub-Riemannian Length Minimizers
by Variational Methods}, preprint 1999.

\bibitem{GP} F.\ Giannoni, P.\ Piccione, {\em An Existence Theory for 
Relativistic Brachistochrones
in Stationary Spacetimes}, J.\  Math.\  Phys.\  {\bf39}, vol.\ 11
(1998), p.\ 6137--6152.

\bibitem{GPV} F.\ Giannoni, P.\ Piccione, J.\ A.\ Verderesi, 
{\em An Approach to the
Relativistic Brachistochrone Problem by sub-Riemannian 
Geometry}, J.\ 
Math.\  Phys.\   {\bf38}, n.\ 12 (1997), 6367--6381.

\bibitem{Kishi} I.\ Kishimoto, {\em The Morse Index Theorem for Carnot--Carath\'eodory
Spaces}, J.\ Math.\ Kyoto Univ.\ {\bf38}, n.\ 2 (1998), 287--293.

\bibitem{L} S.\ Lang, {\em Differential Manifolds}, Springer-Verlag,
Berlin, 1985.

\bibitem{LS} W.\ Liu, H.\ J.\ Sussmann, {\em Shortest Paths for 
sub-Riemannian Metrics of Rank $2$ Distributions}, Memoirs
of the Amer.\ Math.\ Soc.\ {\bf564}, 118 (1995).

\bibitem{M1} R.\ Montgomery, {\em A Survey of Singular Curves in
sub-Riemannian Geometry}, J.\ Dynam.\ Control Systems {\bf1} (1995),
no.\ 1, 49--90.

\bibitem{M2} R.\ Montgomery, {\em Abnormal Minimizers},  SIAM J.\ Control
Optim.\ {\bf32} (1994), no. 6, 1605--1620.

\bibitem{PT1} P.\ Piccione, D.\ V.\ Tausk, {\em On the Banach Differentiable Structure
for Sets of Maps on Non Compact Domains}, to appear in Nonlinear Analysis: Series A,
Theory and Methods.


\end{thebibliography}
\end{document}